# SENSITIVITY ANALYSIS USING PERTURBED-LAW BASED INDICES FOR QUANTILES AND APPLICATION TO AN INDUSTRIAL CASE

ROMAN SUEUR*, BERTRAND IOOSS, THIBAULT DELAGE
*PRISME Department, EDF R&D, 6 quai Watier*
*78401 Chatou CEDEX, France*
email: roman.sueur@edf.fr

In this paper, we present perturbed law-based sensitivity indices and how to adapt them for quantile-oriented sensitivity analysis. We exhibit a simple way to compute these indices in practice using an importance sampling estimator for quantiles. Some useful asymptotic results about this estimator are also provided. Finally, we apply this method to the study of a numerical model which simulates the behaviour of a component in a hydraulic system in case of severe transient solicitations. The sensitivity analysis is used to assess the impact of epistemic uncertainties about some physical parameters on the output of the model.

## 1. Introduction and context of the study

Global sensitivity analysis (GSA) methods are increasingly used tools to understand how the uncertainty in the output of a numerical model is impacted by the uncertainties in the inputs [1]. However, some of the most common of these methods proved inadequate to the study of some output quantities, such as exceedance probabilities, quantiles or even superquantiles (also called expected shortfall or conditional value at risk). This is the case, for example, of the variances based importance measures (Sobol' indices) when applied to the indicator function of a threshold exceedance for the output. To deal with this difficulty and evaluate the behavior of this type of quantity of interest as a function of inputs uncertainties, Lemaître et al. [2] introduced a new family of sensitivity indices called Perturbed-Law based Indices (PLI).
The probability, for a numerical code output, to exceed a threshold is a typical quantity of interest one studies in a structural reliability problem, as it can usually be interpreted as the probability of an undesired event. However, some other quantities, as quantiles or superquantiles, might reveal themselves very useful in the industrial practice. Indeed, a quantile can often be seen as the value associated to a certain level of confidence compared to an expected value, which is often seen as a "best estimate" for the uncertain model output. The associated superquantile then provides a refined measurement of the risk.

In this article, we first present the PLI indices as a general approach for GSA, turning its view to robustness issues. Then we show how these indices can be easily adapted to the study of quantiles, as suggested in [3], and give some useful theoretical results. Finally we present some practical results on a component in a hydraulic system of a power plant submitted to some severe transients which's integrity has to be justified in such a case.

## 2. Perturbed-Law based Sensitivity Indices

We assume to be studying the output $Y$ of a numerical code $G$ which's inputs $\boldsymbol{X} = (X_1, \ldots, X_d) \in \mathbb{R}^d$ are $d$ real-valued independent random variables. We intend to examine the global impact of a change in the values of each $X_i$ to the distribution of $Y$. To this aim, the method developed by Lemaître in [4], consists in perturbing by a $\delta$ quantity a parameter of the marginal density $f_i$ of the $i$-th input of $G$. This perturbation of the $X_i$ input variable induces a change on the output distribution allowing the definition of a sensitivity index which reflects the effect of the distribution of $X_i$ on $Y$.

### 2.1. *PLI for exceedance probabilities*

PLI indices were firstly developed for the study of exceedance probabilities defined by $P = \int_{\boldsymbol{x} \in \mathcal{D}} \mathbb{I}_{\{G(\boldsymbol{x}) > \eta\}} f(\boldsymbol{x}) d\boldsymbol{x}$, where $f$ is the density of the vector $X$, $\mathcal{D}$ is the definition domain of $f$ and $\eta \in \mathbb{R}$ is a threshold value. If $f_{i\delta}$ denotes the perturbed density of $X_i$, $P_{i\delta}$ the subsequent value of the exceedance probability is defined by $P_{i\delta} = \int_{\boldsymbol{x} \in \mathcal{D}} \mathbb{I}_{\{G(\boldsymbol{x}) > \eta\}} f(\boldsymbol{x}) \frac{f_{i\delta}(x_i)}{f_i(x_i)} d\boldsymbol{x}$, where the marginal density of the perturbed input has just been replaced by its perturbed version.

The PLI index $S_{i\delta}$ based on this perturbed probability is defined by:

$$S_{i\delta} = \left(\frac{P_{i\delta}}{P} - 1\right) \mathbb{I}_{\{P_{i\delta} > P\}} + \left(1 - \frac{P}{P_{i\delta}}\right) \mathbb{I}_{\{P_{i\delta} < P\}}. \qquad (1)$$

This definition provides some desirable properties to $S_{i\delta}$, such as:
- $S_{i\delta} = 0$ when $P_{i\delta} = P$ e.g. when perturbing $f_i$ has no impact on the failure probability;
- The sign of $S_{i\delta}$ indicates how the perturbation modifies the failure probability (symetrical behaviour for increase or decrease of the output probability)

### 2.2. *Definition of a density perturbation*

The idea of densities perturbation is to move a specific parameter of the initial density from a $\delta$ quantity. For example, if one choses to introduce a perturbation on the mean value of the density $f_i$, the perturbed density $f_{i,\delta}$ will be such that its mean $\mu_{i,\delta}$ will be equal to $\mu_i + \delta$, where $\mu_i = \mathbb{E}[X_i]$. However, this does not provide an unambiguous definition of $f_{i,\delta}$. So, the definition suggested in [2] and [4] consists in choosing the closest density from $f_i$ in the sense of the Kullback-Leibler divergence, under the constraint that $\mu_{i,\delta} = \mathbb{E}[X_i] + \delta$. More generally, we can give the following formulation for this definition:

$$f_{i\delta} = \underset{\substack{\pi \\ s.t.\ \mathbb{E}_\pi[\psi_k] = \mathbb{E}_{f_i}[\psi_k] + \delta_k \\ k=1,\dots,K}}{\operatorname{argmin}} KL(\pi, f_i)$$

where $\psi_1, \dots, \psi_K$ are $K$ linear constraints on the modified density, and $\delta_1, \dots, \delta_K$ are the values for the perturbations.

### 2.3. *Estimating PLI using Monte-Carlo simulations*

A straightforward approach to calculate PLI values in practice consists in replacing the integral by a discrete sum in the definition of $P$ and $P_{i\delta}$. Hence if we assume having available $\{x^{(n)}\}_{n=1,\dots,N}$ a set of $N$ randomly sampled points in $\mathcal{D}$, independent and identically distributed according to the density $f$, we can define the following Monte-Carlo estimators:

$$\hat{P}^N = \frac{1}{N} \sum_{n=1}^{N} \mathbb{I}_{\{G(x^{(n)})<0\}}, \hat{P}_{i\delta}^N = \frac{1}{N} \sum_{n=1}^{N} L_i^{(n)} \mathbb{I}_{\{G(x^{(n)})<0\}}, \quad (2)$$

where $L_i^{(n)} = f_{i\delta}(x_i^{(n)})/f_i(x_i^{(n)})$ is the likelihood ratio between the perturbed and initial marginal densities of the variable $X_i$. This weighting of each term of the sum is identical to the one used in importance sampling (IS) schemes. For this reason, this method is usually called « reverse importance sampling ». It is noteworthy that the same sample can be used for the calculation of both $\hat{P}^N$ and $\hat{P}_{i\delta}^N$, without need to additional runs of the code $G$.

For a couple $(\hat{P}^N, \hat{P}_{i\delta}^N)$, the plug-in estimator of the $S_{i\delta}$ index is given by :

$$\hat{S}_{i\delta}^N = \left(\frac{\hat{P}_{i\delta}^N}{\hat{P}^N} - 1\right) \mathbb{I}_{\{\hat{P}_{i\delta}^N > \hat{P}^N\}} + \left(1 - \frac{\hat{P}^N}{\hat{P}_{i\delta}^N}\right) \mathbb{I}_{\{\hat{P}_{i\delta}^N < \hat{P}^N\}}.$$

### 3. Adaptation of PLI to the study of the quantiles of the model output

Now, we get interested in another quantity of interest $q^\alpha(Y) = \inf\{t \in \mathbb{R}, F_Y(t) \geq \alpha\}$, the $\alpha$-quantile of the output $Y$, also called value-at-risk (VaR) in the fields of finance or insurance. Provided one is able to define and calculate a perturbed quantile $q_{i,\delta}^\alpha(Y)$, the adaptation of the PLI itself is straightforward. It simply means replacing $P$ and $P_{i\delta}$ in equation (1) by their equivalents $q^\alpha(Y)$ and $q_{i,\delta}^\alpha(Y)$. The main issue is then to be able to estimate these quantiles in a way that guarantees good asymptotic properties to the ensuing quantile-oriented-$\hat{S}_{i\delta}^N$ indice.

### 3.1. *Empirical estimation of a quantile*

The empirical quantile $\hat{q}^{\alpha N}$ is simply the quantity we obtain by replacing the empirical distribution function in the definition of the theoretical quantile: $\hat{q}^{\alpha N} = \inf\{t \in \mathbb{R}, \hat{F}_Y^N(t) \geq \alpha\}$, where $\hat{F}_Y^N(t) = 1/N \sum_{n=1}^{N} \mathbb{I}_{\{y^{(n)} \leq t\}}$. It is worth mentioning that this estimator is strongly consistent and asymptotically normal (see section 2.5 in [5]).

It is obviously possible to apply to this formula the same "reverse importance sampling" principle used for (2). Although the resulting estimator is then given by applying the definition of the quantile to $\tilde{F}_{i\delta}^N = 1/N \sum_{n=1}^{N} L_i^{(n)} \mathbb{I}_{\{y^{(n)} \leq t\}}$, it is preferable to use the following IS estimator of the cumulative distribution function:

$$\hat{F}_{i\delta}^N(t) = \frac{1}{N} \sum_{n=1}^{N} L_i^{(n)} \mathbb{I}_{\{y^{(n)} \leq t\}} \bigg/ \sum_{n=1}^{N} L_i^{(n)},$$

as suggested in [3]. The Figure 1 illustrates how this modifies the initial $F_Y$ estimate in a very simple case of two normal mean-shifted distributions.

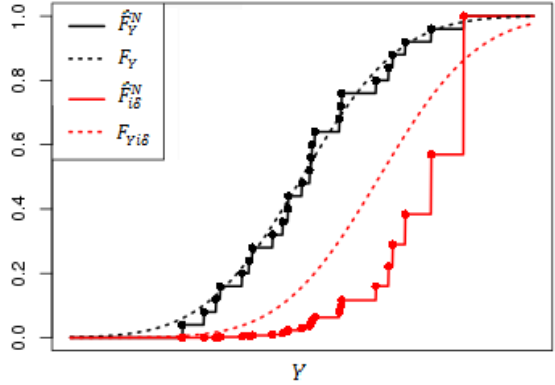

Figure 1. Empirical probability functions $\hat{F}_Y^N$ and $\hat{F}_{i\delta}^N$ and their theoretical models.

The related quantile estimators can be written $\hat{q}_{i\delta}^{\alpha N} = \inf\{t \in \mathbb{R}, \hat{F}_{i\delta}^N(t) \geq \alpha\}$ and $\tilde{q}_{i\delta}^{\alpha N} = \inf\{t \in \mathbb{R}, \tilde{F}_{i\delta}^N(t) \geq \alpha\}$. The quantile-oriented PLI estimator simply becomes:

$$\tilde{S}_{i\delta}^N = \left(\frac{\tilde{q}_{i\delta}^{\alpha N}}{\hat{q}^{\alpha N}} - 1\right) \mathbb{I}_{\{\tilde{q}_{i\delta}^{\alpha N} > \hat{q}^{\alpha N}\}} + \left(1 - \frac{\hat{q}^{\alpha N}}{\tilde{q}_{i\delta}^{\alpha N}}\right) \mathbb{I}_{\{\tilde{q}_{i\delta}^{\alpha N} < \hat{q}^{\alpha N}\}}$$

$$\hat{S}_{i\delta}^N = \left(\frac{\hat{q}_{i\delta}^{\alpha N}}{\hat{q}^{\alpha N}} - 1\right) \mathbb{I}_{\{\hat{q}_{i\delta}^{\alpha N} > \hat{q}^{\alpha N}\}} + \left(1 - \frac{\hat{q}^{\alpha N}}{\hat{q}_{i\delta}^{\alpha N}}\right) \mathbb{I}_{\{\hat{q}_{i\delta}^{\alpha N} < \hat{q}^{\alpha N}\}}.$$

## 3.2. *Asymptotic results about the IS-estimator of the quantile*

The paper [6] provides the demonstration of the a.s. convergence of $\tilde{q}_{i\delta}^{\alpha N}$ to the theoretical quantile as $N \to \infty$, under the assumptions that $Y$ has a positive and differentiable density in a neighborhood of $q^\alpha$, and that there exists a $C \in \mathbb{R}^+$ s.t. $\forall u \in \mathcal{D}_i, f_{i\delta}(u)/f_i(u) < C$, where $\mathcal{D}_i$ is the definition domain of $X_i$. The asymptotic normality of $\tilde{q}_{i\delta}^{\alpha N}$ as $N \to \infty$ is also established:

$$\sqrt{N}\big(\tilde{q}_{i\delta}^{\alpha N} - q^\alpha(Y)\big) \to \frac{\sqrt{Var\left(\mathbb{I}_{\{Y \leq q^\alpha\}} \frac{f_{i\delta}}{f_i}\right)}}{f_Y(q^\alpha)} \mathcal{N}(0,1).$$

This result is also demonstrated in [7] which establishes similar convergences for alternative estimators of $q^\alpha(Y)$, including the already mentioned $\hat{q}_{i\delta}^{\alpha N}$.

## 4. Application to the study of an industrial system's reliability

### 4.1. *Presentation of the application*

The industrial problem we intend to study is the reliability of a component in a production unit which is a part of a hydraulic system. This component, in some incidental situations, might potentially be submitted to a temperature transient leading to a temperature peak (see Figure 2) which could, if too high, involve a risk with regard to the mechanical resistance of the system. Hence, we use a numerical code modeling this kind of transient and giving, as an output quantity, the value of the peak in temperature. However, some physical constants appearing in the physical equations embedded in the code are not perfectly known. These uncertainties, usually called "epistemic uncertainties", induce an uncertainty on the value on the temperature peak.

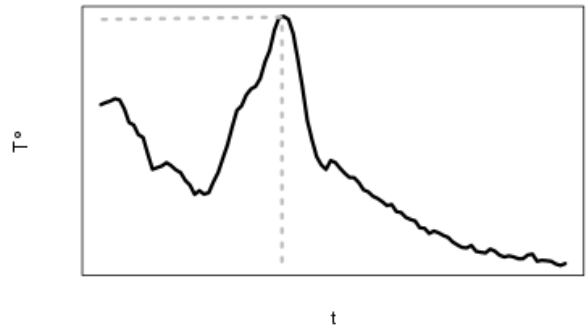

Figure 2. Typical aspect of a temperature transient including a peak

Inputs are not aleatory by nature but are modelled as random variables to represent the lack of knowledge about their true values. For this reason, there is no probabilistic law as such that could objectively represent this uncertainty. There is, for example, no data set which could be used to build a model using common statistical inference methods. The only information available about these physical parameters are bounds obtained by expert judgment. So we choose to use uniform laws between these bounds to perform the uncertainty analysis. However, to examine the robustness of this choice, we first test the impact of a change in mean of these laws. As this code is supposed to be conservative in terms of risk evaluation, we choose to study this impact on the upper-tail 5% quantile of the temperature peak ($\alpha$ =95%).

### 4.2. *Methodology of the study*

As indicated above, we use uniform laws for the nine uncertain inputs of the code. The shape of the perturbed densities for a perturbation on the mean and on the standard deviation (SD) is illustrated in Figure 3. We compute the $\hat{S}_{i\delta_k}^N$ on a regular grid of size $K$ =100 of $\delta_k$ values between -1 and 1, using $N$ =1000 code runs. Note that for a uniform density $U_{[a;b]}$ which's mean is $\mu_U = (a + b)^2/2$ and SD $\sigma_U = |b - a|/\sqrt{12}$, the lower bound of the domain is equal to $a = \mu_U - \sqrt{3}\sigma_U$. To be able to compare the influence of all inputs, we set, for each $X_i$, and each value of $\delta_k$, the perturbed law $f_{i\delta_k}$ as the one defined in section 2.2. such that $\mathbb{E}[f_{i\delta_k}] = \mathbb{E}[f_i] + \delta_k\sqrt{Var(f_i)}$. So it would be theoretically possible to make $\delta_k$ vary from $-\sqrt{3}$ to $+\sqrt{3}$. But when $\mathbb{E}[f_{i\delta_k}]$ is close to the lower bound of the law (which is not changed by the perturbation) numerical problems arise.

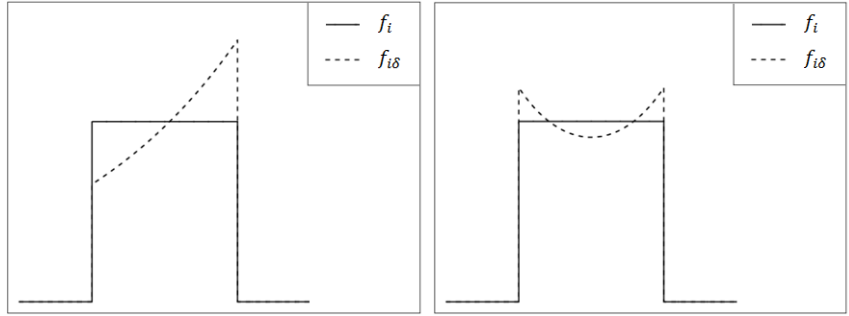

Figure 3. Perturbed density for a uniform initial law and a $\delta$ perturbation on the mean (left) and standard deviation (right).

As the central limit theorem (CLT) given in [7] for $\hat{q}_{i\delta}^{\alpha N}$ requires the estimation of $f_Y(q^\alpha)$ to compute the estimation variance of the IS-quantile estimator, there is no straightforward method to calculates confidence intervals for the $\hat{S}^N_{i\delta_k}$. In particular, the equivalent for quantiles of the CLT for probability oriented-PLI, which is established in [4] (see proposition 3.2.2) cannot be easily proved. To deal with this difficulty, we use a leave one out (LOO) loop over the $N$ simulations to assess confidence intervals (CI).

### 4.3. *Presentation of the results for a perturbation of the mean*

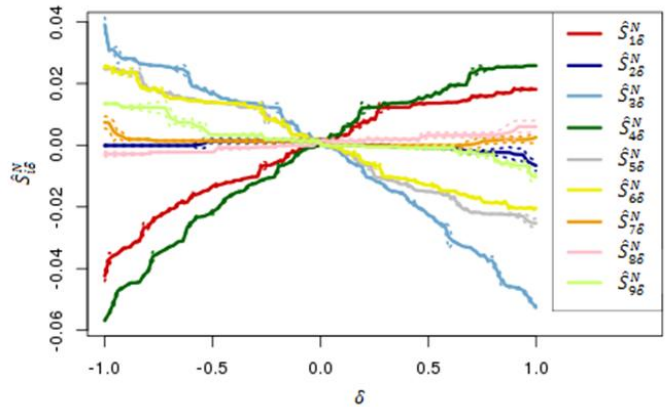

Figure 4. PLI values vs. $\delta$ for the nine inputs of the code

The above Figure 4 shows the values of the PLI estimates (solid lines) with bootstrap CI (dashed lines). The $\alpha$-quantile increases with a positive $\delta$ for variables 1, 4, 7 and 8, and decreases for variables 2, 3, 5, 6 and 9. All indices seem to be monotonic with respect to the perturbation. The most influential inputs on $q^\alpha(Y)$ are variables 1, 4, 5 and 6 and above all variable 3 with a maximum of a 4% increase of the quantile when the mean is diminished of 10% of the SD.

### 5. Perspectives

To overcome the difficulty of estimating the IS-estimator of $q^\alpha(Y)$, and be able to calculate asymptotic CI, a promising track could be to use a kernel estimator for the density of $Y$. Indeed, it would directly provide an estimated value for $f_Y(q^\alpha)$, allowing a plug-in estimation of the variance. In addition to this, the paper [8] gives some useful theoretical results on the convergence of the related quantile estimator. Finally, it would result of a smoothed PLI with respect to $\delta$, which is non continuous when one uses the empirical estimators presented above.